\documentclass[11pt]{article}

\usepackage{setspace}

\usepackage{amsmath}
\usepackage{amssymb}
\usepackage{amsthm}
\usepackage{titlesec}
\usepackage{graphicx}
\usepackage{caption}
\usepackage{subcaption}

\usepackage{diagbox}

\usepackage[T1]{fontenc}
\usepackage[utf8]{inputenc}

\usepackage{comment}

\usepackage{indentfirst}
\usepackage[top=1.25in,left=1.25in,right=1.25in]{geometry}
\newlength{\numone}
\newlength{\widone}
\newlength{\numtwo}
\newlength{\widtwo}
\newtheorem{thm}{Theorem}%[section]
\newtheorem{ex}[thm]{Example}

\numberwithin{equation}{section}

\author{\Large{Riccardo W. Maffucci}}
\newcommand{\Addresses}{{
		\footnotesize
		
		R.W.~Maffucci, \textsc{University of Coventry, United Kingdom CV1}\par\nopagebreak\vspace{-0.35cm}
		\textit{E-mail address}, R.W.~Maffucci: \texttt{riccardowm@hotmail.com}}}

\title{\Large{\uppercase{\bf Rao's Theorem for forcibly planar sequences revisited}}}
\date{}

\def\calA{\mathcal{A}}
\def\calB{\mathcal{B}}
\def\calC{\mathcal{C}}

\def\calP{\mathcal{P}}

\begin{document}
\titleformat{\section}
  {\Large\scshape}{\thesection}{1em}{}
\titleformat{\subsection}
  {\large\scshape}{\thesubsection}{1em}{}
\maketitle
\Addresses

%\pagenumbering{roman}
%\addcontentsline{toc}{section}{Table of contents}
%\tableofcontents

\begin{abstract}
We consider the graph degree sequences such that every realisation is a polyhedron. It turns out that there are exactly eight of them. All of these are unigraphic, in the sense that each is realised by exactly one polyhedron. This is a revisitation of a Theorem of Rao about sequences that are realised by only planar graphs.

Our proof yields additional geometrical insight on this problem. Moreover, our proof is constructive: for each graph degree sequence that is not forcibly polyhedral, we construct a non-polyhedral realisation.
\end{abstract}

\section{Introduction}
This paper deals with finite graphs having no multiple edges or loops. Everywhere we will use the notations $G=(V,E)$ for a graph of order $p$, and
\[s: d_1,\dots,d_p\]
for its degree sequence (i.e., the $d_i$'s are the degrees of the vertices of $G$). The expression $d^k$ in a sequence means $k$ copies of $d$. Vice versa, a sequence of $p$ non-negative integers 
%\[s: d_1,\dots,d_p\]
is called graphical if there exists a graph $G$ with such a degree sequence. Then $G$ is called a realisation of $s$. The Havel-Hakimi algorithm \cite{hakimi,have55} decides whether a given $s$ is graphical, and if so, also constructs a realisation.

We call $s$ \emph{unigraphic} if, up to isomorphism, there is exactly one graph realising $s$. The unigraphic sequences have been characterised in \cite{koren1,li1975}.

Let $\calP$ be a given property  that a graph may or may not satisfy. Examples include, connectivity, planarity, Hamiltonicity, being a line graph, being a perfect graph, \dots. An interesting question is to determine, given a graphical sequence $s$, whether every realisation of it satisfies $\calP$. We then say that $s$ is `forcibly $\calP$'.

The property we will consider is for the graph to be a polyhedron, i.e., the $1$-skeleton of a polyhedral solid. This is exactly the class of $3$-connected, planar graphs \cite{radste}. These may be immersed in the sphere in a unique way. Their regions correspond naturally to polyhedral faces. Two distinct faces are either disjoint, or share exactly one vertex, or share exactly one edge (and its endpoints). Apart from their geometrical and combinatorial interest, they have several applications in the sciences \cite{kuyumc,rouvra,scirih}. A related question to the one considered in this paper, is to characterise sequences that have exactly one polyhedral realisation \cite{delmaf,mafp05}.

Back to the question, which sequences are forcibly polyhedral? It turns out that this problem has just eight solutions.
\begin{thm}
	\label{thm:1}
There are exactly eight forcibly polyhedral sequences, namely
\[s_a: (a+3)^{4-a},(a+2)^a,3^a, \qquad 0\leq a\leq 4,\]
and those corresponding to the octahedron, and square and pentagonal pyramids.
\end{thm}

Our proof of Theorem \ref{thm:1} will be given in section \ref{sec:proof}. We point out that another way to prove Theorem \ref{thm:1} is to use Rao's Theorem \cite{rao978}. In this work, Rao found all forcibly planar sequences. If a sequence is forcibly polyhedral, then it is forcibly planar. It thus remains to establish which solutions listed by Rao are also forcibly $3$-connected. For instance, the sequence of a pyramid is forcibly planar, but as we shall see in section \ref{sec2.1}, all $n$-gonal pyramidal sequences with $n\geq 6$ have realisations that are not $2$-connected.

The forcibly outerplanar sequences were found by Choudum \cite{choudu}. The problem of forcibly $\calP$ sequences has attracted recent attention \cite{baue09,baue13,wang18,wang19,baue20,bar022}.

Our proof of Theorem \ref{thm:1} is independent of the arguments in \cite{rao978,choudu}. We focus on the polyhedral case, and shed some light on the geometric nature of the solutions. For instance, we will see that the polyhedra corresponding to $s_a$ are obtained from the tetrahedron via $a$ `face splitting' operations.

In addition, the nature of our proof is constructive: for each graph degree sequence that is not one of the eight forcibly polyhedral, we find a non-polyhedral realisation.

\section{Proof}
\label{sec:proof}
\paragraph{Intuition.}
Let $G$ be a polyhedral realisation of $s$. We define the new graph
\begin{equation}
\label{eq:tra}
G'=G-v_1v_2-v_3v_4+v_1v_3+v_2v_4,
\end{equation}
for distinct $v_1,v_2,v_3,v_4\in V$. The sequence $s$ is invariant with respect to \eqref{eq:tra}, so that $G'$ also realises $s$. This transformation is one of the main ideas behind the Havel-Hakimi algorithm \cite{hakimi,have55}. Except in a few special cases, we will be able to carefully choose $v_1,v_2,v_3,v_4$ so that $G'$ is not planar, or not $3$-connected.

\subsection{$G$ is not a triangulation}
\paragraph{General strategy.}
Assuming that $G$ is not a triangulation, let $F$ be a face bounded by an $n$-gon,
\[F=[u_1,u_2,\dots,u_n],\]
with $n\geq 4$. The edges $u_1u_2$ and $u_3u_4$ separate $F$ from the faces $F_1$ and $F_2$, say. Consider two vertices $x,y$, with $x\neq u_1,u_2$ on the boundary of $F_1$, and $y\neq u_3,u_4$ on the boundary of $F_2$.

Suppose that $F_1, F_2$ are disjoint faces in $G$. Now $F$ is a face, hence $u_1,\dots, u_n$ is not a separating cycle in $G$. Thereby, there exists in $G$ an $xy$-path $P$ not containing any of $u_1,\dots, u_n$. If $P$ contains more than one vertex from $F_1$ and/or $F_2$, we relabel $x$ to be the last vertex of $F_1$ along $P$, and $y$ the first vertex of $F_2$ along $P$.

\begin{figure}[h!]
	\centering
	\includegraphics[width=4.5cm]{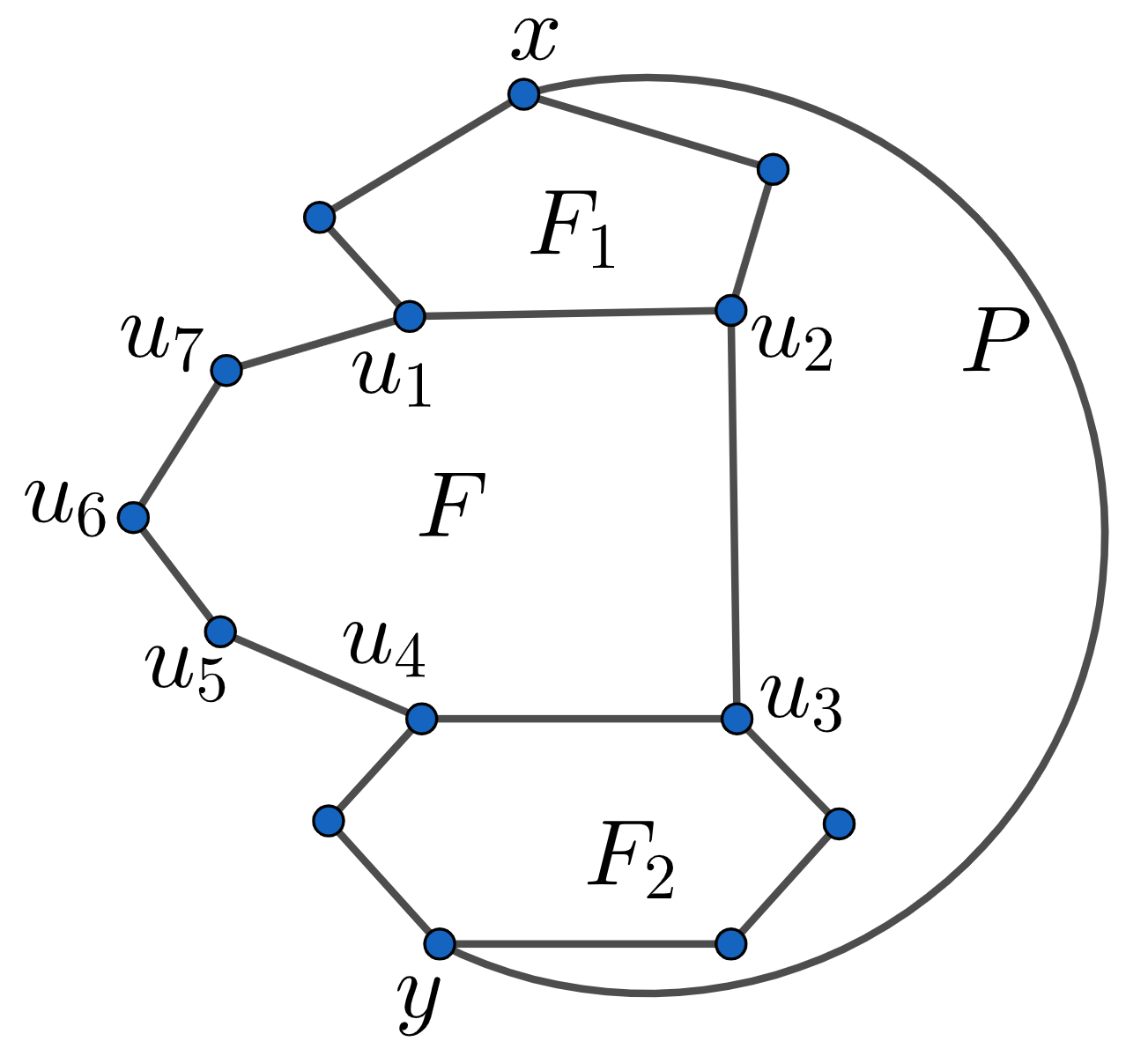}
	\hspace{0.5cm}
	\hspace{0.5cm}
	\includegraphics[width=4.5cm]{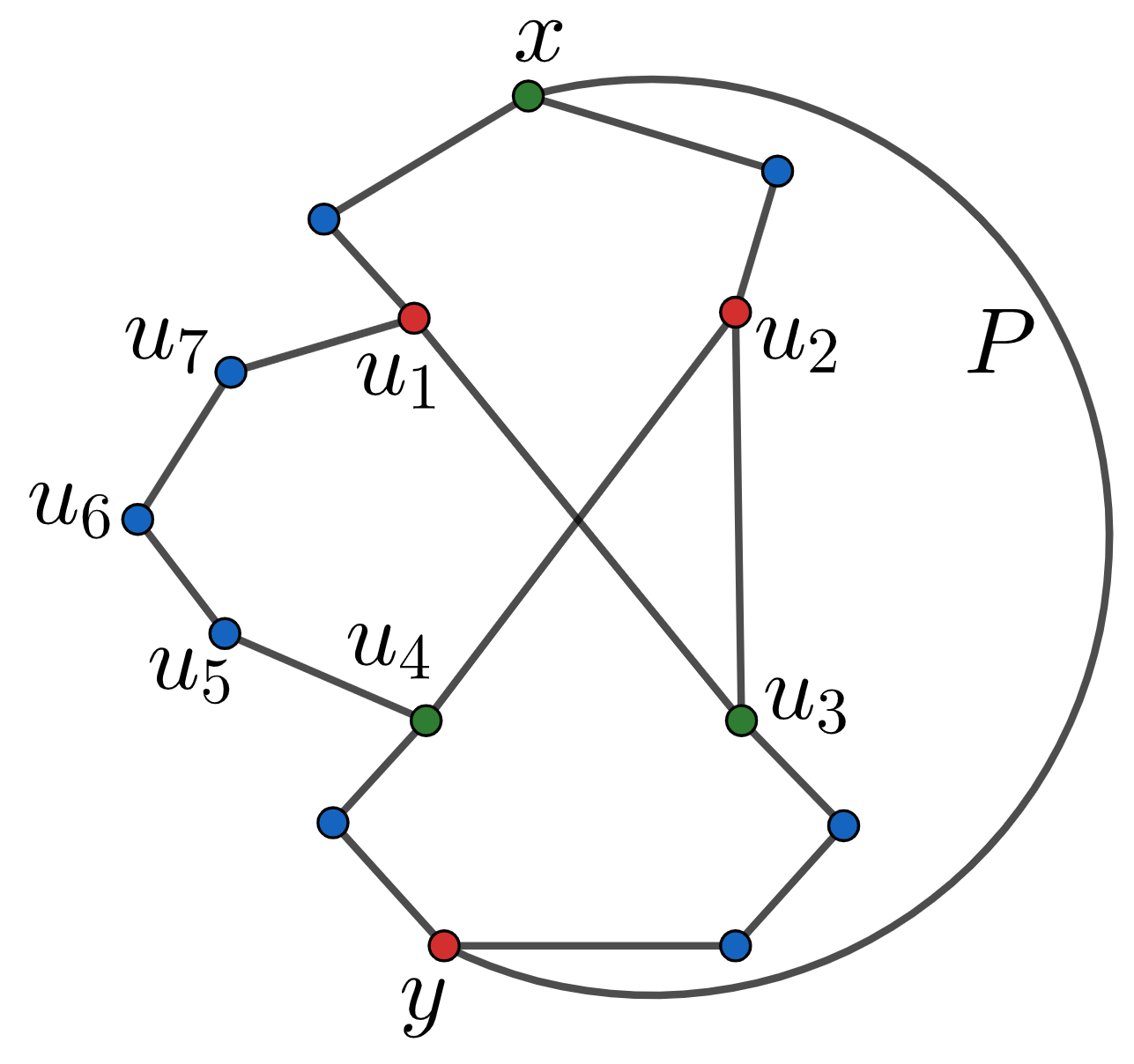}
	\caption{In the generic case, $F_1,F_2$ are disjoint (left). We transform $G$, resulting in the situation on the right. The new graph contains a $K(3,3)$ minor, and is thus non-planar.}
	\label{fig:1}
\end{figure}

We apply to $G$ the transformation \eqref{eq:tra}, with
\[v_i=u_i, \qquad i=1,2,3,4.\]
The resulting $G'$ contains a subgraph homeomorphic to $K(3,3)$, via the partition \[\{\{u_1,u_2,y\},\{u_3,u_4,x\}\}\]
(Figure \ref{fig:1}). Therefore, $s$ is not forcibly polyhedral.

\begin{figure}[h!]
	\centering
	\includegraphics[width=6.5cm]{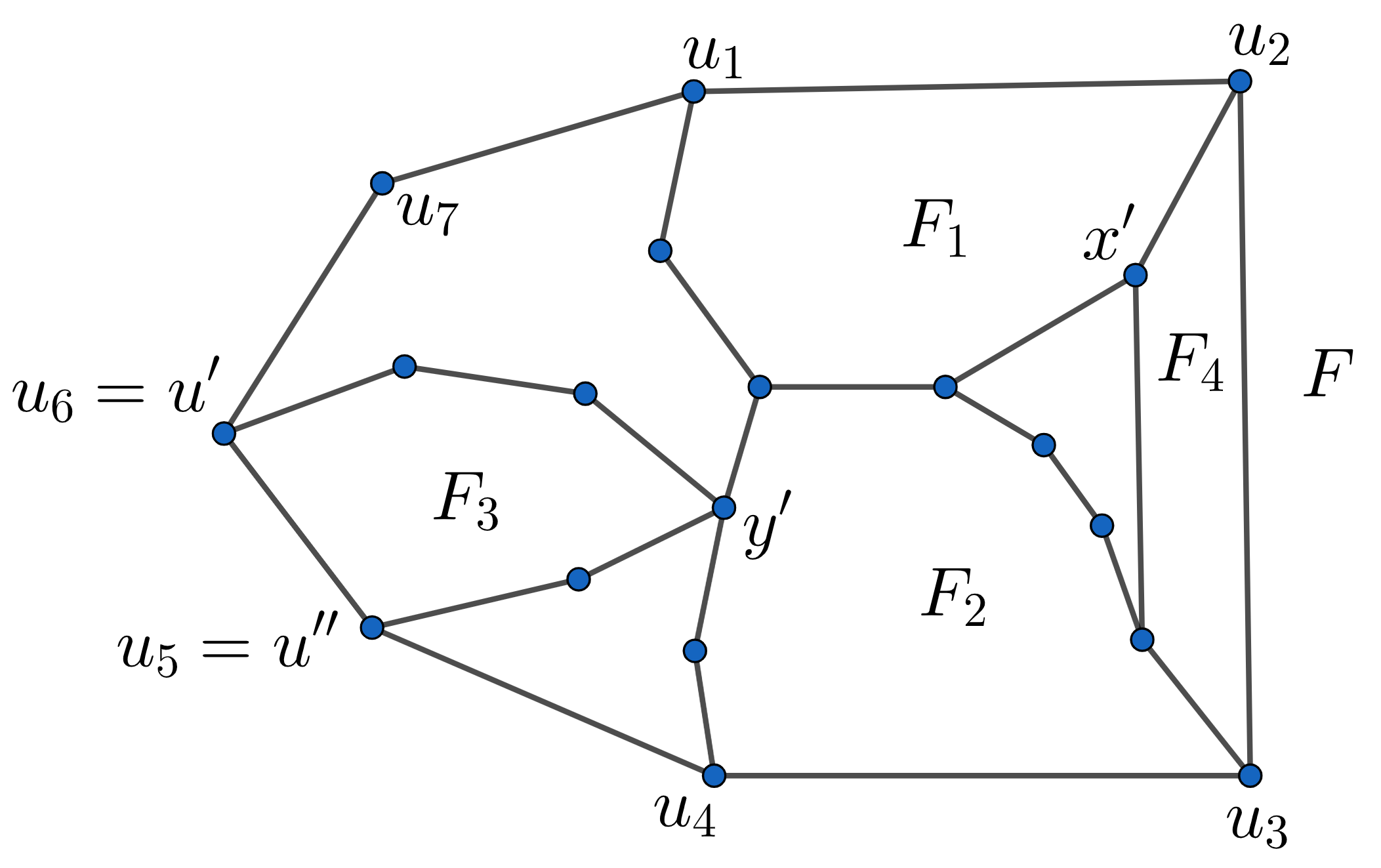}
	\hspace{0.5cm}
	\includegraphics[width=6.5cm]{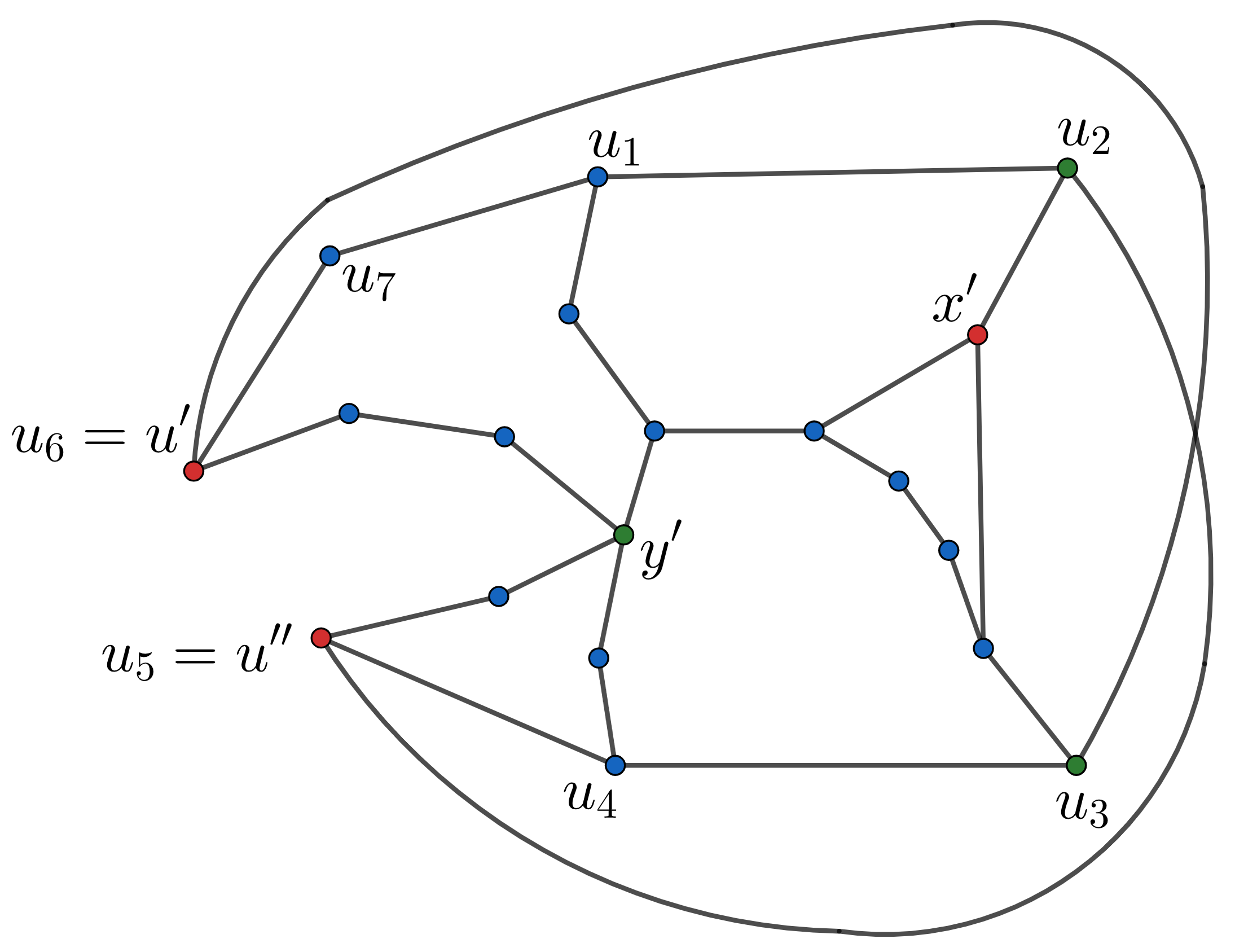}
	\caption{Case of $F_1,F_2$ not disjoint, and $G$ not a pyramid (left). We transform $G$, obtaining a $K(3,3)$ minor in $G'$ (right).}
	\label{fig:2}
\end{figure}

\paragraph{Special cases.}
Now suppose that the faces $F_1, F_2$ are not disjoint. Let $F_4$ be the other face apart from $F$ containing $u_2u_3$. By planarity, unless $G$ is a pyramid, there exists a face $F_3$ adjacent to $F$ and disjoint with $F_4$. Call $u'u''$ the edge separating $F$ from $F_3$ (either $u'u''=u_1u_n$, or $u'u''=u_{m+1}u_{m}$ for some $4\leq m\leq n-1$).
\\
We proceed as above, finding $y'\neq u',u''$ on $F_3$ and $x'\neq u_2,u_3$ on $F_4$ such that there exists in $G$ a $x'y'$-path containing neither other vertices of $F_3,F_4$, nor any vertex of $F$. We then apply \eqref{eq:tra} with
\[v_1=u',\ v_2=u'',\ v_3=u_3,\ v_4=u_2,\]
to obtain a $K(3,3)$ minor in $G'$, with partition
\[\{\{u',u'',x'\},\{u_2,u_3,y'\}\}\]
(Figure \ref{fig:2}). Then $s$ is not forcibly polyhedral.

If $G$ is an $n$-gonal pyramid, $n\geq 6$, with central vertex $c$ and base given by the cycle $u_1,u_2,\dots,u_n$, we take \eqref{eq:tra} with
\[v_1=u_2,\ v_2=u_3,\ v_3=u_6,\ v_4=u_5.\]
The resulting graph $G'$ is not $2$-connected, as $G'-c$ is disconnected.

On the other hand, the square and pentagonal pyramids are the only graphs realising their respective sequences, thus $4,3^4$ and $5,3^5$ are forcibly polyhedral.

\begin{ex}
Applying the above algorithm to the triangular prism, we get the graph $K(3,3)$.
\end{ex}

\label{sec2.1}

\subsection{$G$ is a triangulation}
\paragraph{General strategy.}
Similarly to section \ref{sec2.1}, we inspect a general scenario, leaving out special cases for later. Suppose that there exist two vertices $u,w$ of $G$ of degree at least four, that are not adjacent. We denote by $U,W$ their respective neighbourhoods. Since $G$ is a triangulation, we may write
\[U=\{u_1,\dots,u_n\},\]
where
\begin{equation}
	\label{eqn:cycle}
u_1,\dots,u_n
\end{equation}
is a cycle (an analogous statement may be made for the vertex $w$ of course).

Next, by $3$-connectivity, there are at least three internally disjoint $wu$-paths in $G$. These are chosen so that their combined length is minimal, and will be denoted by $\calA,\calB,\calC$. Each contains exactly one element of $U$, as penultimate vertex. We call $a,b,c\in U$ the three penultimate vertices, and $a_1,b_1,c_1$ the neighbours of $w$, along $\calA,\calB,\calC$ respectively. Note that possibly one or more of $a=a_1$, $b=b_1$, $c=c_1$ hold.%and $d$ an element of $U$ distinct from $a,b,c$.

\begin{figure}[h!]
	\centering
	\includegraphics[width=3.5cm]{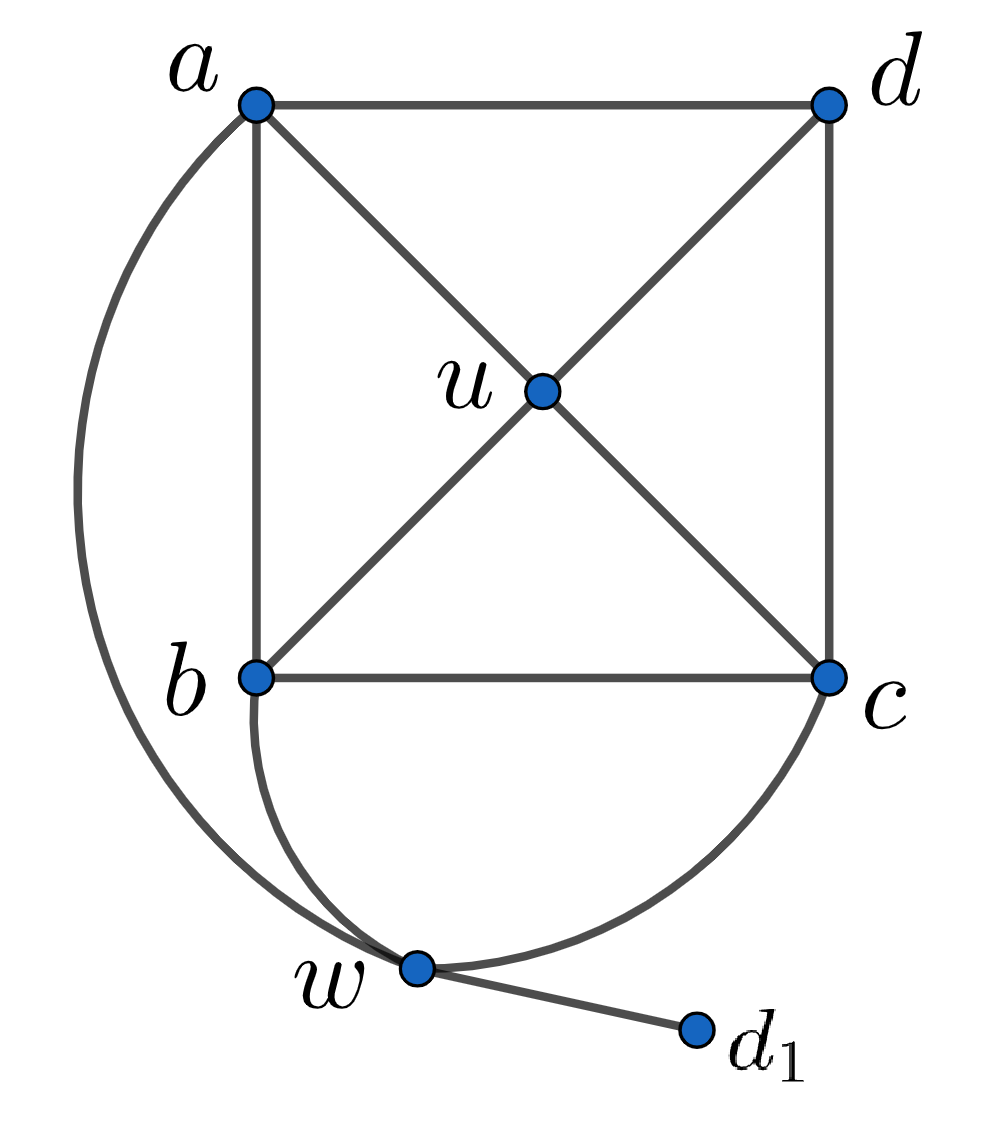}
	\hspace{1.0cm}
	\includegraphics[width=3.5cm]{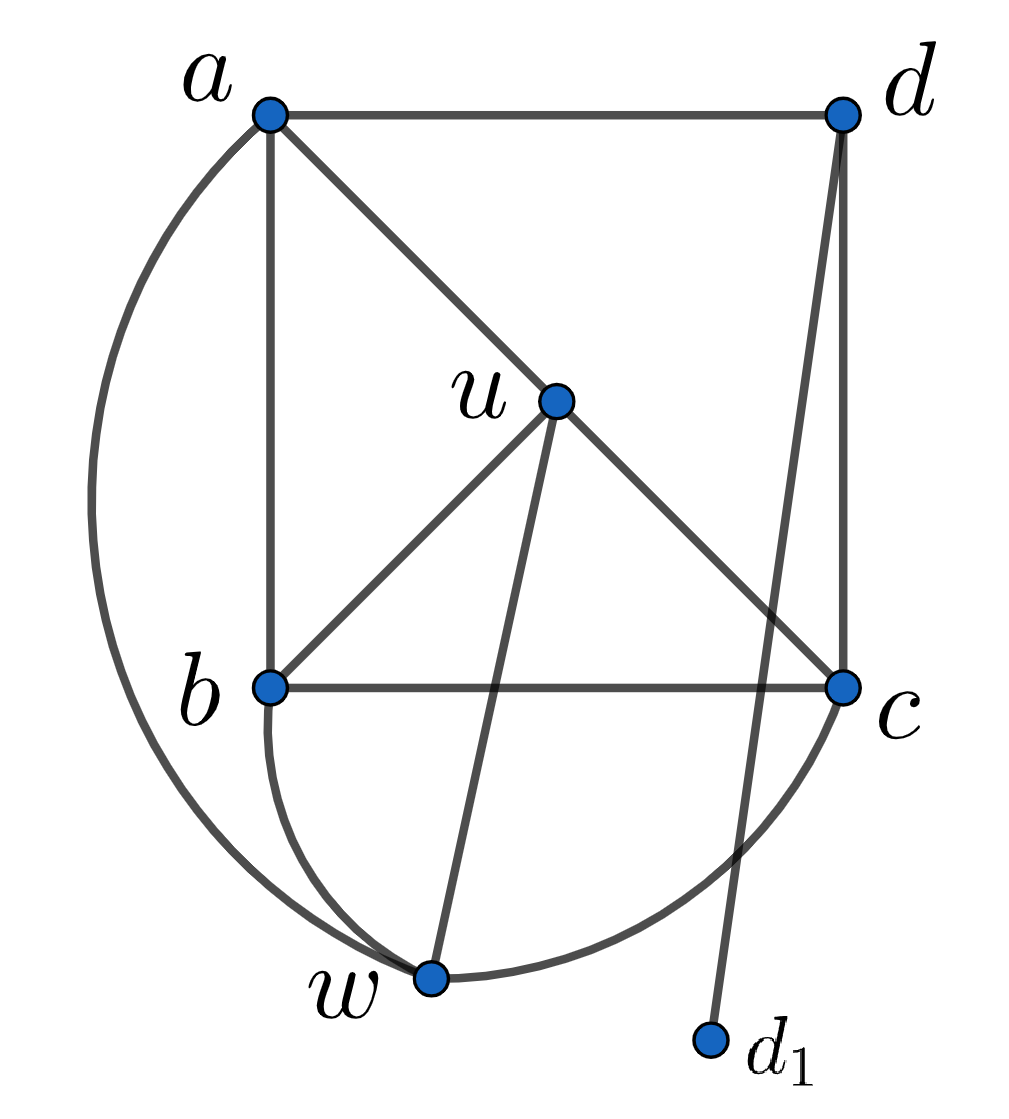}
	\caption{The typical strategy for $G$ maximal planar, and $u,w$ vertices of degree at least four that are not adjacent.}
	\label{fig:3}
\end{figure}

Since $u,w$ are of degree at least four, there exist
\[d\in U\setminus\{a,b,c\}, \qquad d_1\in W\setminus\{a_1,b_1,c_1\}.\]
Suppose for the moment that
\[
d\neq d_1 \ \text{ and } \ dd_1\not\in E(G).
\]
We define $G'$ as in \eqref{eq:tra}, with
\[v_1=u,\ v_2=d,\ v_3=w,\ v_4=d_1.\]
Then the vertices
\[\{a,b,c,u,w\}\]
determine in $G'$ a subgraph homeomorphic from $K_5$, thus $G'$ is not planar, and in particular $s$ is not forcibly polyhedral (Figure \ref{fig:3}).

Now assume instead that
\[
d=d_1 \ \text{ or } \ dd_1\in E(G).
\]
We claim that in this case there are at least four internally disjoint $wu$-paths in $G$. Indeed, $\calA,\calB,\calC$ cannot contain $d$, as each includes exactly one element of $U$. None of them may contain $d_1$ either, as this would contradict the minimality of their combined length. By the way, for the same reason, $\calA,\calB,\calC$ are actually each of length $2$ or $3$. W.l.o.g., $a,b,c,d$ appear in this order around the cycle \eqref{eqn:cycle}. If $a\neq a_1$, say, then we define
$G'$ as in \eqref{eq:tra}, with
\[v_1=a,\ v_2=a_1,\ v_3=c,\ v_4=x,\]
where
\begin{equation*}
x=\begin{cases}
c_1 & c\neq c_1\\w & c=c_1.
\end{cases}
\end{equation*}
Again we obtain a copy of $K_5$ in $G'$, this time via the vertices
\[\{a,b,c,d,u\},\]
hence $s$ is not forcibly polyhedral.

The only other possibility is that $a=a_1$, $b=b_1$, $c=c_1$, and $d=d_1$ all hold. Moreover, by the same argument, any further vertices adjacent to $u$ must be adjacent to $w$ as well (and vice versa). Now assume that $G$ contains a vertex $y$ not adjacent to $u,w$. Then by planarity, there exist two elements of $U$ such that each $yu$-path contains at least one of them. This contradicts $3$-connectivity. It follows that $G$ is the $n$-gonal bipyramid, with base given by the cycle \eqref{eqn:cycle}.

Finally, if $n\geq 5$, we take \eqref{eq:tra} with
\[v_1=u,\ v_2=u_1,\ v_3=w,\ v_4=u_3.\]
Then
\[\{u_2,u_4,u_5,u,w\}\]
forms a copy of $K_5$ in $G'$, thus $s$ is not forcibly polyhedral. On the other hand, the sequence $4^6$ is unigraphic. Its realisation is the octahedron.

\paragraph{Special cases.}
It remains to establish what happens if in a triangulation $G$ all vertices of degree at least four are pairwise adjacent. Let
\[m=m(G):=|\{v\in V: \deg(v)\geq 4\}|.\]
Firstly, we deduce that $m\leq 4$, otherwise $G$ would contain a copy of $K_5$. If $m=0$, then $G$ is $3$-regular in addition to being a triangulation. In other words, $G$ is simply the tetrahedron. If $m=1$, then $s: n,3^{n}$, i.e. $G$ is the $n$-gonal pyramid, $n\geq 4$, impossible. For $m=2$, we write
\[s: \deg(x),\deg(y),3^{p-2}.\]
Now
\[\deg(x)+\deg(y)+3(p-2)=6p-12.\]
On the other hand, $\deg(x),\deg(y)\leq p-1$, implying $3p-6\leq 2(p-1)$, i.e. $p\leq 4$, impossible.

Similarly, if $m=3$, we obtain that $G$ has $p$ vertices, of which three have degree $p-1$, and the rest degree $3$. Due to planarity, the only candidate solution is the sequence $s_1: 4,4,4,3,3$ of the triangular bipyramid. It is unigraphic, hence forcibly polyhedral.

Lastly, let $m=4$, i.e.
\[s: \deg(x),\deg(y),\deg(z),\deg(w),3^{p-4}.\]
On one hand, $x,y,z,w$ are pairwise adjacent, hence the total number of edges of type $uv$ with $u\in\{x,y,z,w\}$ and $v$ of degree $3$ is given by
\[\deg(x)+\deg(y)+\deg(z)+\deg(w)-12=(6p-12)-3(p-4)-12=3p-12.\]
On the other hand, the sum of degrees of the vertices of degree $3$ in $G$ is exactly
\[3(p-4)=3p-12.\]
Combining these considerations, each vertex of degree three in $G$ is adjacent to exactly three of $x,y,z,w$. Since $w$ is adjacent to all of $x,y,z$, for $G$ not to contain a $K(3,3)$ minor, there can be at most one vertex of degree three adjacent to all of $x,y,z$, and similarly for the other three triples. In particular, it must hold that $p\leq 8$. Now the symmetries of the tetrahedron generated by $\{x,y,z,w\}$ tell us that the case $m=4$ has exactly three candidate solutions, namely
\[s_2: 5,5,4,4,3,3, \quad s_3: 6,5,5,5,3,3,3, \quad s_4: 6,6,6,6,3,3,3,3.\]
In other words, for $0\leq a\leq 4$, the polyhedron corresponding to $s_a$ is precisely the tetrahedron with $a$ `face splittings' (the splitting of a polyhedral face is an operation defined by inserting a new vertex inside the face, together with edges to all vertices of the face). One checks that $s_2,s_3,s_4$ are all unigraphic, and the proof of Theorem \ref{thm:1} is complete.

\bibliographystyle{abbrv}
\bibliography{bibgra}

\begin{thebibliography}{10}

\bibitem{bar022}
A.~Bar-Noy, T.~B{\"o}hnlein, D.~Peleg, and D.~Rawitz.
\newblock On realizing a single degree sequence by a bipartite graph.
\newblock In {\em 18th Scandinavian Symposium and Workshops on Algorithm Theory
  (SWAT 2022)}. Schloss Dagstuhl-Leibniz-Zentrum f{\"u}r Informatik, 2022.

\bibitem{baue13}
D.~Bauer, H.~J. Broersma, J.~van~den Heuvel, N.~Kahl, and E.~Schmeichel.
\newblock Toughness and vertex degrees.
\newblock {\em Journal of graph theory}, 72(2):209--219, 2013.

\bibitem{baue09}
D.~Bauer, S.~L. Hakimi, N.~Kahl, and E.~Schmeichel.
\newblock Sufficient degree conditions for k-edge-connectedness of a graph.
\newblock {\em Networks: An International Journal}, 54(2):95--98, 2009.

\bibitem{baue20}
D.~Bauer, L.~Lesniak, A.~Nevo, and E.~Schmeichel.
\newblock On the necessity of {C}hv{\'a}tal’s {H}amiltonian degree condition.
\newblock {\em AKCE International Journal of Graphs and Combinatorics},
  17(2):665--669, 2020.

\bibitem{choudu}
S.~Choudum.
\newblock Characterization of forcibly outerplanar graphic sequences.
\newblock In {\em Combinatorics and Graph Theory: Proceedings of the Symposium
  Held at the Indian Statistical Institute, Calcutta, February 25--29, 1980},
  pages 203--211. Springer, 1981.

\bibitem{delmaf}
J.~Delitroz and R.~W. Maffucci.
\newblock On unigraphic polyhedra with one vertex of degree $ p-2$.
\newblock {\em arXiv preprint arXiv:2301.08021}, 2023.

\bibitem{hakimi}
S.~Hakimi.
\newblock On the realizability of a set of integers as degrees of the vertices
  of a graph.
\newblock {\em SIAM Journal Applied Mathematics}, 1962.

\bibitem{have55}
V.~Havel.
\newblock A remark on the existence of finite graphs.
\newblock {\em Casopis Pest. Mat.}, 80:477--480, 1955.

\bibitem{koren1}
M.~Koren.
\newblock Sequences with a unique realization by simple graphs.
\newblock {\em Journal of Combinatorial Theory, Series B}, 21(3):235--244,
  1976.

\bibitem{kuyumc}
A.~Kuyumcu and A.~Garcia-Diaz.
\newblock A polyhedral graph theory approach to revenue management in the
  airline industry.
\newblock {\em Computers \& industrial engineering}, 38(3):375--395, 2000.

\bibitem{li1975}
S.-Y.~R. Li.
\newblock Graphic sequences with unique realization.
\newblock {\em Journal of Combinatorial Theory, Series B}, 19(1):42--68, 1975.

\bibitem{mafp05}
R.~W. Maffucci.
\newblock Characterising $3$-polytopes of radius one with unique realisation.
\newblock {\em arXiv preprint arXiv:2207.02725}, 2022.

\bibitem{rao978}
S.~Rao.
\newblock Characterization of forcibly planar degree sequences.
\newblock {\em ISI Tech. Report}, (36/78), 1978.

\bibitem{rouvra}
D.~Rouvray.
\newblock Graph theory in chemistry.
\newblock {\em Royal Institute of Chemistry, Reviews}, 4(2):173--195, 1971.

\bibitem{scirih}
I.~Sciriha and P.~W. Fowler.
\newblock Nonbonding orbitals in fullerenes: Nuts and cores in singular
  polyhedral graphs.
\newblock {\em Journal of chemical information and modeling}, 47(5):1763--1775,
  2007.

\bibitem{radste}
E.~Steinitz and H.~Rademacher.
\newblock Vorlesungen {\"u}ber die {T}heorie der {P}olyeder. {S}pringer 1934.

\bibitem{wang18}
K.~Wang.
\newblock An efficient algorithm to test forcibly-connectedness of graphical
  degree sequences.
\newblock {\em arXiv preprint arXiv:1803.00673}, 2018.

\bibitem{wang19}
K.~Wang.
\newblock Forcibly-biconnected graphical degree sequences: Decision algorithms
  and enumerative results.
\newblock {\em Theory and Applications of Graphs}, 6(2):4, 2019.

\end{thebibliography}

\end{document}